\documentclass[a4paper,11pt]{article}
\usepackage{amssymb}

\newcommand{\R}{\mathbb{R}}
\newcommand{\Z}{\mathbb{Z}}
\newcommand{\N}{\mathbb{N}}

\newcommand{\beq}{\begin{equation} }
\newcommand{\eqq}{\end{equation} }
\newcommand{\cuad}{{\sqcap\kern-.68em\sqcup}}

\newtheorem{teo}{Theorem}[section]

\newtheorem{proposition}{Proposition}[section]
\newtheorem{example}{Example}[section]

\newtheorem{lemma}{Lemma}[section]
\newtheorem{corollary}{Corollary}[section]
\newtheorem{remark}{Remark}[section]
\newcommand{\bremark}{\begin{remark} \em}
\newcommand{\eremark}{\end{remark} }

\def\beeq{\begin{equation}}
\def\eeq{\end{equation}}
\newcommand{\begeqaet}{\begin{eqnarray*}}
\newcommand{\eneqaet}{\end{eqnarray*}}

\begin{document}

\begin{center}{\bf  \Large   Sharp embedding of Sobolev spaces involving general kernels and its application }\medskip
\bigskip\medskip

{\bf Huyuan Chen\footnote{chenhuyuan@yeah.net} \quad \quad  Hichem Hajaiej\footnote{hichem.hajaiej@gmail.com}}
\smallskip

$^1$Department of Mathematics, Jiangxi Normal University,
 Nanchang, Jiangxi 330022, PR China\\
and\\
 $^2$Department of Mathematics, College of Science, King Saud University P.O. Box 2455,
 Riyadh 11451, Saudi Arabia.

\bigskip
\begin{abstract}
The purpose of this paper is to extend the embedding theorem of  Sobolev spaces involving general kernels and  we provide
a sharp critical exponent in these embeddings. As an application, solutions for equations driven by a general integro-differential operator, with homogeneous Dirichlet boundary conditions, is established by
using the Mountain Pass Theorem.

\end{abstract}
\end{center}
\medskip

\vspace{1mm}

\medskip

\date{}

\setcounter{equation}{0}
\section{ Introduction}

In the study of weak solutions for semilinear elliptic equations,
 the embedding from corresponding Sobolev space to $L^q$ space
 plays a fundamental role, especially the compact embedding. In a recent work,  Di Nazza,  Palatucci and  Valdinoci in \cite{EGE} made
  a clear description for the fractional Sobolev space $W^{s,p}(\Omega)$ and gave an elegant proof for
 the embedding theorem from $W^{s,p}(\Omega)$ to $L^q(\Omega)$, which is continuous when $q\in[1,\frac{Np}{N-sp}]$ and
  compact when $q\in[1,\frac{Np}{N-sp})$, where $s\in(0,1)$, $sp<N$ and $\Omega\subset \R^N$ is a bounded  extension domain  with $N\geq2$.

 Motivated by the above work, our purpose of this paper is to build a sharp embedding theorem of Sobolev space involving
 general kernel $K$ and by using this embedding theorem to search for  weak solutions  to problem
\begin{equation}\label{RE}
 \arraycolsep=1pt
\begin{array}{lll}
 \mathcal{L}_Ku+f(x,u)=0\quad & \rm{in}\quad\Omega,\\[2mm]
 \phantom{    \mathcal{L}_K+u^p--}
u=0\quad & \rm{in}\quad \Omega^c,
\end{array}
\end{equation}
where  $\Omega\subset \R^N$ is an open bounded  $C^2$ domain  with $N\geq2$ and the nonlocal operator $\mathcal{L}_K$ is defined by
$$\mathcal{L}_Ku(x)=\frac12\int_{\R^N}[u(x+y)+u(x-y)-2u(x)]K(y)dy$$
with the kernel $K:\R^N\setminus\{0\}\to (0,+\infty)$  satisfying
\begin{equation}\label{1}
\int_{\R^N}\min\{|x|^2,1\}K(x)dx<+\infty
\end{equation}
and
\begin{equation}\label{3}
K(x)=K(-x),\qquad x\in \R^N\setminus\{0\}.
\end{equation}
Moreover, we assume that $K$ is decreasing monotone in the following sense
\begin{equation}\label{monotone}
K(x)\geq K(y)\qquad {\rm{if}}\ \ |x|\leq|y|.
\end{equation}
A typical example for $K$ is given by $K(x) = |x|^{-(N+2s)}$ with $s\in(0,1)$ and then $\mathcal{L}_K$ is the fractional
Laplacian operator $-(-\Delta)^s$.

During the last  years, non-linear equations involving  general
integro-differential operators, especially, fractional Laplacian,
have been studied by many authors.
Caffarelli and Silvestre \cite{CS} studied the fractional Laplacian through
extension theory. Caffarelli and Silvestre  \cite{CS1,CS2}, Ros-Oton and Serra  \cite{RS} investigated regularity results  for fractional elliptic
equations. Sire and Valdinoci in \cite{S1},  Felmer and Wang in \cite{FW}, Hajaiej \cite{HH1,HH2} and Felmer, Quaas and Tan \cite{FQT} obtained symmetry property of
solutions for  semilinear equation involving the fractional Laplacin.
More interests on fractional elliptic equations see \cite{CFQ,CV2,CV3,CV4,HH} and the references therein.

Recently, Servadei and  Valdinoci in \cite{RE} obtained a solution of (\ref{RE}) via  Mountain Pass Theorem
 under the hypothesis that there exist $\lambda>0$ and
$s\in (0,1)$ such that
$$
K(x)\geq \lambda |x|^{-(N+2s)},\quad x\in\R^N\setminus\{0\}
$$
and nonlinear term $f$
is subcritical, that is, $$|f(x,t)|\leq a_1+a_2|t|^{q-1}\quad {\rm a.e.}\ x\in\Omega,\ t\in\R $$ with
$q\in(2,\frac{2N}{N-2s})$ and constants $a_1,a_2>0$. We say that $\frac{2N}{N-2s}$ is the critical exponent, denoted by $2^*(s)$.

In this paper, we are also interested in studying problem (\ref{RE}) with more general kernels and our purpose is to find new criterion for critical exponent, where  we could deal with the following case
\begin{equation}\label{OS22}
\liminf_{|x|\to0^+}K(x)|x|^{N}\in(0,\infty).
\end{equation}
 To this end, we define
\begin{equation}\label{case 1}
s_0=\sup\{s\ge0\ |\ \lim_{r\to0^+}r^{2s}\int_{B_r^c(0)}K(y)dy=+\infty\}.
\end{equation}
We remark that if $K$ satisfies (\ref{1}) and is nonnegative, then
the definition in (\ref{case 1}) is equivalent to
$$s_0=\sup\{s\ge0\ |\ \lim_{r\to0^+}r^{2s}\int_{B_1(0)\setminus B_r(0)}K(x)dx=+\infty\}$$
By the fact that $\int_{B_1^c(0)}K(x)dx$ is bounded from (\ref{1}).

Our basic setting is that $s_0>0.$
In section 2, we will prove that $s_0\le 1$ and exhibit  an  example in which the kernel $K$ satisfying (\ref{OS22}) makes $s_0\in(0,1)$.
We note that the limit of $r^{2s_0}\int_{B_r^c(0)}K(y)dy$, as $r\to 0$, could be in $[0, \infty]$ or even no exists.
Denote
\begin{equation}\label{case 1.1}
l_\infty:=  \liminf_{r\to0^+}r^{2s_0}\int_{B_r^c(0)}K(y)dy,
\end{equation}
then it occurs one of the cases:  Case 1: $l_\infty=0$ and Case 2: $l_\infty\in(0,\infty]$,

 Our first aim is to study the Sobolev space involving general kernel $K$. Denote by $X$  the linear
space of Lebesgue measurable functions from $\R^N$ to $\R$ such that
the restriction to $\Omega$ of any function $g$ in $X$ belongs to $L^2(\Omega)$ and
$$\int_{\R^{2N}\setminus\mathcal{O}}(g(x)-g(y))^2K(x-y)dxdy<+\infty,$$
where $\mathcal{O}:=\Omega^c\times\Omega^c$. The space $X$ is  endowed
with the norm as
\begin{equation}\label{Norm}
\|g\|_X=(\|g\|^2_{L^2(\Omega)}+\int_{\R^{2N}\setminus\mathcal{O}}(g(x)-g(y))^2K(x-y)dxdy)^{1/2}.
\end{equation}
Now we define the following Sobolev space
  $$X_0=\{g\in X\ | \ g=0\ \ {\rm a.e.\ in}\
\Omega^c\}$$
equipped the norm (\ref{Norm}). From (\ref{1}), we stress that $C^2_0(\Omega)\subseteq X_0,$
see \cite{RE}, and so $X$ and $X_0$ are nonempty.

Now we are ready for an embedding theorem.

\begin{teo}\label{th 1}
Assume that $K$ satisfies (\ref{1}), (\ref{monotone}), (\ref{case
1}) with  $s_0\in(0,1]$, $2^*(s_0)=\frac{2N}{N-2s_0}$ and $l_\infty$ is defined by (\ref{case 1.1}).
 Then $(X_0, \|\cdot\|_{X})$ is a Hilbert space and \\[1mm]
 $(i)$ if $l_\infty=0$,  the embedding
\begin{equation}\label{embedding 01}
X_0\hookrightarrow L^{q}(\Omega)
\end{equation} is continuous and compact for $q\in[1,2^*(s_0))$. Moreover, for $q\in[1,2^*(s_0))$
there exists $C>0$ such that
\begin{equation}\label{sobolev inequality}
\|g\|_{L^{q}}\le C\|g\|_X,\quad \forall g\in X_0;
\end{equation}
$(ii)$ if $l_\infty\in(0,\infty]$,   the embedding
(\ref{embedding 01}) is continuous for $q\in[1,2^*(s_0)]$ and compact for $q\in[1,2^*(s_0))$, and
the embedding inequality (\ref{sobolev inequality}) holds for for $q\in[1,2^*(s_0)]$.
\end{teo}

\begin{example}\label{ex 0}
Let
\begin{equation}\label{1.01}
K(x)=\frac1{|x|^{N+2s_0}}\left[(-\log|x|)_++1\right]^\sigma,\quad x\in\R^N\setminus\{0\},
\end{equation}
where $\sigma\in\R$ and $(-\log|x|)_+=\max\{-\log|x|,0\}$.
When $s_0\in(0,1)$   $\sigma\in\R$ or $s_0=1$ $\sigma<-1$, the kernel $K$ defined by (\ref{1.01}) satisfies (\ref{1}) and (\ref{monotone}).

We note that $l_\infty=0$ if $\sigma<0$, $l_\infty\in(0,\infty)$ if $\sigma=0$ and
$l_\infty=\infty$ if $\sigma>0$.
In particular, $s_0\in(0,1)$ and $\sigma=0$,   the embedding (\ref{embedding 01}) coincides the results in \cite{EGE}.
Especially, when $s_0=1$ and $\sigma<-1$, $2^*(s_0)=2^*$ the critical exponent for $H^1_0(\Omega)\Subset L^{2^*}(\Omega)$.
\end{example}

Now we are able to make use of  Theorem \ref{th 1} to study the existence of weak
solutions of (\ref{RE}).
Before stating the existence result we make precise the definition of weak solution that we use in the article.
We say that a function $u\in X_0$ is a weak solution of (\ref{RE}) if
\begin{equation}\label{WRE}
 \int_{\R^N\times\R^N}[u(x)-u(y)][\varphi(x)-\varphi(y)]K(x-y)dxdy=\int_\Omega f(x,u(x))\varphi(x) dx,
\end{equation}
for any $\varphi\in X_0$.

The existence result can be  stated as follows.

\begin{teo}\label{th 2}
Assume that $f(x,u)=|u|^{p-2}u$,  $K$ satisfies (\ref{1}-\ref{monotone}), $2^*(s_0)=\frac{2N}{N-2s_0}$, where
$s_0\in(0,1]$ defined in (\ref{case
1}).

Then  problem (\ref{RE}) admits a nontrivial  weak solution  for $p\in(2,2^*(s_0))$.
\end{teo}

\begin{remark}
Take $K$ as example \ref{ex 0} with $s_0\in(0,1)$ and $\sigma\in\R$ or $s_0=1$ and $\sigma<-1$, then problem (\ref{RE}) admits a weak solution for $f(x,u)=|u|^{p-2}u$ with $p\in(2,2^*(s_0))$.

Take $K$ as example \ref{ex 1},  problem (\ref{RE}) admits a weak solution for $f(x,u)=|u|^{p-2}u$ with $p\in(2,2^*(s_0))$.

\end{remark}

The paper is organized as follows. In Section 2, we
analyze some basic properties of the kernel $K$ and give an example showing that
$s_0$ makes sense. In Section 3, we study the Sobolev embedding theorem in our setting.
Finally, we prove the existence of weak solution to (\ref{RE}) in Section 4.

\setcounter{equation}{0}
\section  {Discussion to the kernel $K$}

This section is devoted to  the properties of the kernel $K$.

\begin{proposition}\label{prop 2}
Assume that  $s_0$ is defined by (\ref{case 1}) and $K$ satisfies (\ref{1}). Then
$(i)$ for any $s<s_0$, we have $$\lim_{r\to0^+} r^{2s}\int_{B_r^c(0)}K(y)dy=+\infty;$$
$(ii)$
\begin{equation}\label{7}
s_0\le \inf\{s\ge 0\ | \ \lim_{r\to0^+} r^{2s}\int_{B_r^c(0)}K(y)dy=0\}\le 1;
\end{equation}
$(iii)$  if there exists $s_1\le s_2$ such that
\begin{equation}\label{1.33}
\liminf_{|x|\to0^+}K(x)|x|^{N+2s_1}>0\quad{\rm and}\quad \limsup_{|x|\to0^+}K(x)|x|^{N+2s_2}<\infty,
\end{equation}
 then  $s_0\in[s_1,s_2]$.
\end{proposition}
{\bf Proof.}
\emph{$(i)$}  By the definition of $s_0$, there at least are a sequence of positive numbers $\{s_n\}$ such that
$$s_n<s_0,\quad \lim_{n\to\infty}s_n=s_0,\quad \lim_{r\to0^+}r^{2s_n}\int_{B_r(0)}K(y)dy=+\infty.$$
Then for any $s<s_0$, there exists $s_n$ such that $s<s_n$ and then
$$\lim_{r\to0^+}r^{2s}\int_{B_r(0)}K(y)dy\ge\lim_{r\to0^+}r^{2s_n}\int_{B_r(0)}K(y)dy=+\infty.$$

\emph{$(ii)$}  By (\ref{1}) and $K$ being nonnegative, we have that for any $r\in(0,1)$,
\begin{eqnarray*}
  \infty&>&\int_{\R^N}\min\{|x|^2,1\}K(x)dx \\
  &>&\int_{B_1(0)\setminus B_r(0)} |x|^2K(x)dx+ \int_{\R^N\setminus B_1(0)}K(x)dx\\
   &\ge&  r^2\int_{\R^N\setminus B_r(0)}K(x)dx.
\end{eqnarray*}
Then for any $s>1$, we have that
\begin{eqnarray*}
r^{2s}\int_{B_r^c(0)}K(x)dx=r^{2(s-1)}[r^2\int_{B_r^c(0)}K(x)dx]\to0\quad {\rm as}\ r\to0.
\end{eqnarray*}
Thus, $\inf\{s\ge 0\ | \ \lim_{r\to0^+} r^{2s}\int_{B_r^c(0)}K(y)dy=0\}\le 1$.\smallskip

We now prove the first inequality (\ref{7}). We denote
$$s_{00}=\inf\{s\ge 0\ | \ \lim_{r\to0^+} r^{2s}\int_{B_r^c(0)}K(y)dy=0\}.$$
 Since for any
$s> s_{00}$, we have that
$$\lim_{r\to0^+} r^{2s}\int_{B_r^c(0)}K(x)dx=0.$$
By the definition of $s_0$, we have $s_0\le s$ and then by arbitrary of $s>s_{00}$, we obtain that $s_0\le s_{00}$.\smallskip

\emph{$(iii)$}   For any $s<s_1$ by (\ref{1.33}), we have
\begin{eqnarray*}
r^{2s}\int_{B_r^c(0)}K(x)dx&=&r^{2(s-s_1)}[r^{2s_1}\int_{B_r^c(0)}K(x)dx]
\\&\ge&r^{2(s-s_1)}\inf_{|x|\in(0,1)}(K(x)|x|^{N+2s_1})[r^{2s_1}\int_r^1\tau^{-2s_1-1}d\tau]
\\&\ge&r^{2(s-s_1)}\int_r^1\tau^{-1}d\tau\inf_{|x|\in(0,1)}(K(x)|x|^{N+2s_1})
\\&\to&\infty\quad {\rm as}\ r\to0.
\end{eqnarray*}

By the definition of $s_0$, we have $s_0\ge s$ and then by arbitrary of $s<s_1$, we obtain that $s_0\ge s_1$.
Similarly to prove $s_0\le s_2$. \hfill$\Box$

\begin{lemma}\label{lm 5}
$(i)$  Assume that the kernel $K$ satisfies  (\ref{monotone}) and  is
 continuous in $\R^N\setminus\{0\}$, then
$K$ is radially symmetric about the origin.

$(ii)$  Assume that the kernel $K$ satisfies (\ref{1}), (\ref{monotone}) and (\ref{case
1}) with $s_0>0$. Then for any $s\in(0,s_0)$,
 there exists a sequence $\{r_n\}$ of positive numbers which converges to 0  and
\begin{equation}\label{3.2}
  \lim_{r_n\to0^+}r_n^{N+2s}\inf_{|x|=r_n}K(x)=+\infty.
\end{equation}
\end{lemma}
{\bf Proof.}
$(i)$ By contradiction, we may assume that there exist
$x_1,y_1\in \R^N\setminus\{0\}$ such that $|x_1|=|y_1|$ and $K(x_1)>K(y_1)$.
Since $K$ is continuous in $\R^N\setminus\{0\}$, then there exists $x_2\in \R^N\setminus\{0\}$ such that $|x_2|>|x_1|$ and
$$K(x_2)\ge K(x_1)-\frac{K(x_1)-K(y_1)}{2}>K(y_1),$$ which is impossible with the assumption (\ref{monotone}). \smallskip

$(ii)$ By Proposition \ref{prop 2} $(i)$, we have that for $s\in(0,s_0)$ and $\epsilon\in(0,s_0-s)$,
\begin{equation}\label{2.1}
  \lim_{r\to0^+}r^{2(s+\epsilon)}\int_{B_r^c(0)}K(x)dx=+\infty.
\end{equation}

 Let $\tilde K(r)=\inf_{|x|=r}K(x)$, then by (\ref{monotone}), we have $\tilde K(r_1)\le \tilde K(r_2)$ for $r_1\ge r_2$ and
$K(x)\le \tilde K(r)$ for any $|x|>r$.

 If (\ref{3.2}) doesn't hold, then there no exist any sequence $\{r_n\}$ converging to zero such that (\ref{3.2}) holds, that is
   $$\limsup_{r\to0^+}r^{N+2s}\tilde K(r)<+\infty.$$
 Together with $\tilde K$ is decreasing, then there exists $C>0$ such that
 $$\tilde K(r)\le Cr^{-N-2s},\quad r\in(0,1).$$
For any $x\in B_1(0)\setminus\{0\}$, we have $K(x)\le \tilde K(\frac{|x|}2)$,
\begin{eqnarray*}
  r^{2(s+\epsilon)}\int_{B_1(0)\setminus B_r(0)}K(x)dx &\le &  r^{2(s+\epsilon)}\int_{B_1(0)\setminus B_r(0)}\tilde K(\frac{|x|}2)dx\\
  &\le& C2^{N+2s}r^{2(s+\epsilon)} \int_r^1\tau^{-1-2s}d\tau
  \\&\le& Cr^{2\epsilon}.
\end{eqnarray*}
Together with (\ref{1}), we have
$$
 \lim_{r\to0^+}r^{2(s+\epsilon)}\int_{B_r^c(0)}K(x)dx=0.
$$
which contradicts with (\ref{2.1}). The proof is complete.
\hfill$\Box$

\medskip

To end this section, we construct an example of $K$ satisfying (\ref{OS22}) for which $s_0\in(0,1)$.\\[1mm]
\begin{example}\label{ex 1}  Let
\begin{equation}\label{je 3}
K(x)=\left\{ \arraycolsep=1pt
\begin{array}{lll}
 a_n^{-N-2s},\ \ |x|\in[a_{n+1},a_n),\\[2mm]
|x|^{-N},\ \ |x|\in[a_1,1),\\[2mm]
|x|^{-N-2s},\ \ |x|\in[1,+\infty).
\end{array}
\right.
\end{equation}
where $s\in(0,1)$, $a_0\in (0,1)$, $a_n=a_0^{b^n}$ with $n\in \mathbb{N}$ and $b=\frac{N+2s}{N}$.\\[1mm]
Then $$\liminf_{r\to0^+}K(r)r^{N}=1\quad{\rm and}\quad s_0\in(0,s).$$
\end{example}
 {\bf Proof.}
We  observe that $\lim_{n\to+\infty}a_n=0$ and
\begin{eqnarray*}
  K(a_n)a_n^{N} = a_{n-1}^{-N-2s} a_n^{N}
   = a_0^{-b^{n-1}(N+2s)} a_n^{N}
   = (a_0^{b^{n}})^{-N}a_n^{N}=1,
\end{eqnarray*}
then we have
$$\liminf_{r\to0^+}K(r)r^{N}=1.$$
Combining Proposition \ref{prop 2} $(iii)$ and the fact of $\limsup_{r\to0^+}K(r)r^{N+2s}\le1$, we have that
$$s_0\in[0,s).$$

 Now we prove that  $s_0>0$.
For $r\in(0,a_1)$, there exists $n\in \mathbb{N}$
such that $a_{n+1}\le r< a_n.$ If $n$ big enough, we have $a_{n}\le \frac12 a_{n-1}$.
Combining with $b>1$, then
\begin{eqnarray*}
\int_{B_{a_1(0)}\setminus B_r(0)}K(y)dy&=&|\omega_N|[
(a_n-r)^Na_n^{-N-2s}+\sum_{k=2}^{n}(a_{k-1}-a_k)^Na_{k-1}^{-N-2s}]
\\&\geq&|\omega_N|\sum_{k=2}^{n}(a_{k-1}-a_k)^Na_{k-1}^{-N-2s}
\\&\ge&|\omega_N|2^{-N} a_{n-1}^{-2s},
\end{eqnarray*}
where $w_N$ is the unit sphere of $\R^N$.
Choose $\beta=b^{-2}s>0$,  then we obtain that
$$a_{n-1}^{-2s}\ge a_{n+1}^{-2\beta}.$$
Therefore,
$$\liminf_{r\to0^+}r^{2\beta}\int_{B_{a_1(0)}\setminus B_r(0)}K(y)dy\ge 2^{-N}|\omega_N|.$$
By Proposition \ref{prop 2} $(iii)$, we obtain that $s_0\ge \beta>0.$ \hfill$\Box$

\setcounter{equation}{0}
\section  {Sobolev spaces }

 In this section, we will consider some embedding results inspired from \cite{EGE}. First we
 introduce some basic spaces and   some useful tools to  prove embedding theorems.

\begin{lemma}\label{lm 1}
 Assume that $K$ satisfies (\ref{1}), (\ref{monotone}),  (\ref{case 1}) with $s_0>0$ and $l_\infty$ is defined by (\ref{case 1.1}).
Let $x\in \R^N$ and  $E\subset\R^N$ be a measurable set with
$|E|\in(0, +\infty)$, then\\
$(i)$ if $l_\infty=0$, for any $s\in(0,s_0)$, there exists $C>0$ such that
\begin{equation}\label{3.1}
  \int_{E^c}K(x-y)dy\geq C|E|^{-\frac{2s}N};
\end{equation}
$(ii)$ if $l_\infty\in(0,\infty]$,
  there exists $C>0$ such that (\ref{3.1}) holds with $s\in(0,s_0]$.
\end{lemma}
{\bf Proof.} We just need to prove that the conclusion of Lemma \ref{lm 1}
holds for a sequence $E_n$ with $|E_n|>0$ and
$\lim_{n\to\infty}|E_n|=0$. Let
$\rho_n=(\frac{|E_n|}{\omega_N})^{1/N}$, then it follows that
$|E_n^c\cap B_{\rho_n}(x)|=|E_n\cap B_{\rho_n}^c(x)|$. Therefore,
by (\ref{monotone}),
we have that
$$K(x-y)\ge \inf_{|z|=\rho_n}K(z), \quad y\in E_n^c\cap B_{\rho_n}(x),$$
$$K(x-y)\le \inf_{|z|=\rho_n}K(z), \quad y\in E_n\cap \bar B_{\rho_n}^c(x).$$
Thus
\begin{eqnarray}
 \int_{E_n^c}K(x-y)dy &=& \int_{E_n^c\cap B^c_{\rho_n}(x)}K(x-y)dy+\int_{E_n^c\cap B_{\rho_n}(x)}K(x-y)dy \nonumber\\
   &\ge&  \int_{E_n^c\cap B^c_{\rho_n}(x)}K(x-y)dy+\inf_{|z|=\rho_n}K(z)|E_n^c\cap B_{\rho_n}(x)|\nonumber\\
   &\ge&  \int_{E_n^c\cap B^c_{\rho_n}(x)}K(x-y)dy+\inf_{|z|=\rho_n}K(z)|E_n\cap \bar B_{\rho_n}^c(x)|\nonumber\\
   &=&   \int_{B^c_{\rho_n}}K(x-y)dy.\label{20-09-0}
\end{eqnarray}
\emph{ $(i)$} By Proposition \ref{prop 2} $(i)$ and $s_0>0$, we observe that for any $s\in(0,s_0)$
\begin{equation}\label{20-09-1}
\lim_{r\to0^+}r^{2s}\int_{B_r^c(0)}K(y)dy=\infty.
\end{equation}
Then by (\ref{20-09-0}), there exists $C>0$ such that
$$ \int_{E_n^c}K(x-y)dy\ge C|E_n|^{-\frac{2s}N}.$$

\noindent\emph{ $(ii)$} Since $l_\infty>0$, then there exists $\sigma\in(0,1)$ such that for $r\in(0,1)$
$$r^{2s_0}\int_{B_r^c(0)}K(y)dy\ge \sigma l_\infty,$$ which, together with (\ref{20-09-0}), implies  that
$$ \int_{E_n^c}K(x-y)dy\ge \sigma l_\infty|E_n|^{-\frac{2s_0}N}.$$
For $s\in(0,s_0)$, it is the same as the proof of $(i)$.
   \hfill $\Box$

\begin{lemma}\label{lm 2} \cite[Lemma 6.2]{EGE}
Assume that $s\in(0,1)$, $2s<N$ and $T>1$. Let $n\in \Z$ and  $\{a_k\}$
be a bounded, nonnegative, decreasing sequence  with $a_k=0$ for any
$k\geq n$.
Then,
\begin{eqnarray*}
\sum_{k\in\Z}a_k^{1-\frac{2s}{N}}T^k\leq C\sum_{k\in \Z, a_k\not=0}a_{k+1}a_k^{-\frac{2s}{N}}T^k,
\end{eqnarray*}
for a suitable constant $C=C(s,T,N)>0$, independent of $n$.
\end{lemma}

\begin{lemma}\label{lm 3}
Assume that $K$ satisfies (\ref{1}), (\ref{monotone}),  (\ref{case 1}) with $s_0\in(0,1)$ and $l_\infty$ is defined by (\ref{case 1.1}).
 Let  $f\in L^\infty(\R^N)$ be compactly
supported, then $$\int_{\R^{2N}}|f(x)-f(y)|^2K(x-y)dxdy\geq
C\sum_{k\in \Z,
a_k\not=0}a_{k+1}a_k^{-\frac{2s}N}2^{2k},$$
where $a_k=|\{|f|>2^k\}|$, $k\in \Z$, $C=C(N,K)>0$ and  the choice of $s$ is the same as in Lemma \ref{lm 1}.
\end{lemma}
{\bf Proof.} The proof is similar to Lemma 6.3  in \cite{EGE} just replaced the kernel by $K$. For reader's convenience, we give the detail below.  Firstly, we assume that $f$ is nonnegative. If not, we
replace $f$ by $|f|$. Let $A_k:=\{f>2^k\}$, $D_k:=A_k\setminus
A_{k+1}$, $d_k:=|D_k|$ and $S:= \sum_{j\in \Z,
a_{j-1}\not=0}2^{2j}a_{j-1}^{-\frac{2s}N}d_j.$
Then
\begin{equation}\label{4}
\{(i,j)\in\Z^2\ s. t.\ a_{i-1}\not=0\ {\rm and}\
a_{j-1}^{-\frac{2s}N}d_j\not=0\}\subset \{(i,j)\in\Z^2\ s. t.\
a_{j-1}\not=0\}.
\end{equation}
Then we have that
\begin{eqnarray*}
\sum_{i\in \Z, a_{i-1}\not=0}\sum_{j\in \Z,
j\ge
i+1}2^{2i}a_{i-1}^{-\frac{2s}N}d_j&=&
\sum_{i\in \Z, a_{i-1}\not=0}\sum_{j\in \Z, j\ge
i+1,a_{i-1}^{s}d_j\not=0}2^{2i}a_{i-1}^{-\frac{2s}N}d_j\\&\leq&\sum_{i\in
\Z}\sum_{j\in \Z, j\ge
i+1,a_{i-1}\not=0}2^{2i}a_{i-1}^{-\frac{2s}N}d_j
\\&=&\sum_{j\in
\Z,a_{j-1}\not=0}\sum_{i\in \Z, i\le
j-1}2^{2i}a_{i-1}^{-\frac{2s}N}d_j
\\&\leq&\sum_{j\in
\Z,a_{j-1}\not=0}\sum_{i\in \Z, i\le
j-1}2^{2i}a_{j-1}^{-\frac{2s}N}d_j
\\&=&\sum_{j\in
\Z,a_{j-1}\not=0}\sum_{k=0}^{+\infty}2^{2j-2}2^{-2k}a_{j-1}^{-\frac{2s_0}N}d_j
\\&\leq& S.
\end{eqnarray*}
Fixed $i\in\Z$ and $x\in D_i$, for any $l\in\Z$ with
$l\leq i-2$ and any $y\in D_l$, we have that $$|f(x)-f(y)|\geq
2^{i-1}$$ and therefore,
\begin{eqnarray*}
\sum_{l\in\Z,l\leq i-2}\int_{D_j}|f(x)-f(y)|^2K(x-y)dy&\geq&
2^{2i-2}\sum_{l\in\Z,l\leq
i-2}\int_{D_j}K(x-y)dy\\&=&2^{2i-2}\int_{A^c_{i-1}}K(x-y)dy.
\end{eqnarray*}
By Lemma \ref{lm 1}, we have $$\sum_{l\in\Z,l\leq
i-2}\int_{D_l}|f(x)-f(y)|^2K(x-y)dy\geq c_0
2^{2i}a_{i-1}^{-\frac{2s}N},$$ for some
suitable $c_0>0$. As a consequence, for any $i\in\Z$,
$$\sum_{l\in\Z,l\leq
i-2}\int_{D_i\times D_l}|f(x)-f(y)|^2K(x-y)dxdy\geq c_0
2^{2i}a_{i-1}^{-\frac{2s}N}d_i$$ and then,
\begin{eqnarray*}
\sum_{i\in\Z,a_{i-1}\not=0}\sum_{l\in\Z,l\leq i-2}\int_{D_i\times
D_l}|f(x)-f(y)|^2K(x-y)dxdy\geq c_0 S.
\end{eqnarray*}
Thus, we obtain
\begin{eqnarray*}
&&\sum_{i\in\Z,a_{i-1}\not=0}\sum_{l\in\Z,l\leq i-2}\int_{D_i\times
D_l}|f(x)-f(y)|^2K(x-y)dxdy\\&&\geq c_0
[\sum_{i\in\Z,a_{i-1}\not=0}2^{2i}a_{i-1}^{-\frac{2s}N}a_i-\sum_{i\Z,a_{i-1}\not=0}\sum_{j\in\Z,j\ge
i+1}2^{2i}a_{i-1}^{-\frac{2s}N}d_j]\\&&\geq
c_0(2^{2i}a_{i-1}^{-\frac{2s}N}a_i-S).
\end{eqnarray*}
So, it follows that
\begin{eqnarray*}
&&\int_{\R^N\times
\R^N}|f(x)-f(y)|^2K(x-y)dxdy\\&\geq&2\sum_{i\in\Z,a_{i-1}\not=0}\sum_{l\in\Z,l\leq
i-2}\int_{D_i\times D_l}|f(x)-f(y)|^2K(x-y)dxdy\\&\geq&
c_0(\sum_{i\in\Z,a_{i-1}\not=0}2^{2i}a_{i-1}^{-\frac{2s}N}a_i).
\end{eqnarray*}
\hfill $\Box$

\begin{lemma}\label{lm 4}
Assume that $q\in [1,+\infty)$,  $f:\R^N\to \R$ is a measurable function.
For any $n\in\N$,
$$f_n(x):=\max\{\min\{f(x),n\},-n\},\quad \forall x\in\R^N.$$
Then $$\lim_{n\to+\infty}\|f_n\|_{L^q(\R^N)}=\|f\|_{L^q(\R^N)}.$$
\end{lemma}
{\bf Proof. } The details of the proof refers to \cite[Lemma 6.4]{EGE} or \cite{BH}.\hfill$\Box$

Now we can give the statement of embedding theorem as follows:

\begin{teo}\label{sobolev}
 Assume that $K$ satisfies (\ref{1}), (\ref{monotone}),  (\ref{case 1}) with $s_0>0$ and $l_\infty$ is defined by (\ref{case 1.1}).
 Then\\[1mm]
 $(i)$ if $l_\infty=0$, then for $s\in(0,s_0)$ there exists
$C>0$ such that for any  $f\in X_0$, we have
\begin{equation}\label{imbeding}
\|f\|_{L^{2^*(s)}(\Omega)}\leq C(\int_{\R^N}\int_{\R^N}|f(x)-f(y)|^2K(x-y)dxdy)^{\frac12};
\end{equation}
 $(ii)$ if $l_\infty\in(0,\infty]$,
then   (\ref{imbeding}) holds with $s=s_0$.
\end{teo}
{\bf Proof.} First we note that
\begin{equation}\label{5}
\int_{\R^N}\int_{\R^N}|f(x)-f(y)|^2K(x-y)dxdy<+\infty.
\end{equation}
Without loss of generality, we can assume that $f\in
L^\infty(\R^N)$. Indeed, let $f_n$ be defined as in  Lemma \ref{lm 4}, then
combining with Lemma \ref{lm 4} and  (\ref{5}), we make use of  the
Dominated Convergence Theorem to imply
$$\lim_{n\to\infty}\int_{\R^{2N}}|f_n(x)-f_n(y)|^2K(x-y)dxdy=\int_{\R^{2N}}|f(x)-f(y)|^2K(x-y)dxdy,$$
which allows us to obtain estimate for function $ f\in X_0$.

Take $s$, $a_k$ and $A_k$ defined as in Lemma \ref{lm 3}, then we
have that
 \begin{eqnarray*}
\|f\|^{2^*(s)}_{L^{2^*(s)}(\R^N)}=\sum_{k\in\Z}\int_{A_k\setminus
A_{k+1}}|f(x)|^{2^*(s)}dx\leq \sum_{k\in \Z}2^{2^*(s)(k+1)}a_k,
\end{eqnarray*}
that is,
$$\|f\|^{2}_{L^{2^*(s)}(\R^N)}\leq 4(\sum_{k\in \Z}2^{2^*(s)k}a_k)^{2/2^*(s)}.$$
Since $2<2^*(s)$, then
$$\|f\|^{2}_{L^{2^*(s)}(\R^N)}\leq 4\sum_{k\in \Z}2^{2k}a_k^{2/2^*(s)}.$$
By  Lemma \ref{lm 2} with $T=4$, it follows that
$$\|f\|^{2}_{L^{2^*(s)}(\R^N)}\leq C\sum_{k\in \Z}2^{2k}a_{k+1}a_k^{-\frac{2s}{N}}.$$
for a suitable constant $C$ depending on $N,K$.

Finally, it suffices to apply Lemma \ref{lm 4} to obtain that
$$\|f\|_{L^{2^*(s)}(\R^N)}\leq C(\int_{\R^N}\int_{\R^N}|f(x)-f(y)|^2K(x-y)dxdy)^{\frac12},$$
up to relabeling the constant $C$. Since $f\in X_0$, $f=0$ in $\Omega^c$, then (\ref{imbeding}) holds.
 \hfill $\Box$

\begin{corollary}\label{co 1}
The norm (\ref{Norm}) in $X_0$ is equivalent to
\begin{equation}\label{norm 1}
  \|f\|_{X_0}:=(\int_{\R^N}\int_{\R^N}|f(x)-f(y)|^2K(x-y)dxdy)^{\frac12}.
\end{equation}

\end{corollary}
{\bf Proof.} We only need to prove that there exists $C>0$ such that for any $f\in X_0$,
$$\|f\|_X\le C\|f\|_{X_0}.$$
It follows by Theorem \ref{sobolev} that
\begin{eqnarray*}
&&\|f\|_X^2 = \int_\Omega f^2(x)dx+\int_{\R^N}\int_{\R^N}|f(x)-f(y)|^2K(x-y)dxdy\\
 &&\le  |\Omega|^{1-\frac2{2^*(s)}}(\int_\Omega |f|^{2^*(s)}(x)dx)^{\frac2{2^*(s)}}+\int_{\R^N}\int_{\R^N}|f(x)-f(y)|^2K(x-y)dxdy
 \\&&\le C\int_{\R^N}\int_{\R^N}|f(x)-f(y)|^2K(x-y)dxdy.
 \end{eqnarray*}
The proof is complete.\hfill$\Box$

\begin{teo}\label{compact imbedding}
Assume that $K$ satisfies (\ref{1}), (\ref{monotone}),  (\ref{case 1}) with $s_0>0$ and
$\mathcal{T}$ is a bounded subset of $X_0$.

Then $\mathcal{T}$ is pre-compact in $L^{q}(\Omega)$, $q\in[1,2^*(s_0))$.
\end{teo}
\textbf{Proof.}
We first prove that $\mathcal{T}$ is pre-compact in $L^2(\Omega)$. To this end, we only show that $\mathcal{T}$ is totally bounded in $L^2(\Omega)$.
By Lemma \ref{lm 5}$(ii)$, there exists $\{r_n\}$ positive and convergent to 0
such that
$$
\lim_{n\to\infty}r_n^NK(r_n)=+\infty.
$$
  Let $\rho:\R_+\to \{\frac{r_n}2,n\in\N\}$ such that, denoting $\rho_\epsilon=\rho(\epsilon)$, for any $\epsilon>0$, $\rho_\epsilon=r_n$ for some $n$ and
  $$\lim_{\epsilon\to0^+} \rho_\epsilon=0.$$ It is obvious that
\begin{equation}\label{21-09-0}
\lim_{\epsilon\to0^+}(2\rho_\epsilon)^NK(2\rho_\epsilon)=+\infty.
\end{equation}

 Let
$ \eta_\epsilon=\epsilon \rho_\epsilon^{ \frac N2}$ and take a collection of disjoints cubes $Q_1,....,Q_M$ of side $\rho_\epsilon$
such that
$$\Omega\subset \bigcup_{j=1}^N Q_j.$$
For any $x\in\Omega$, there exists a unique integer $j(x)$ in $\{1,...,M\}$ such that $x\in Q_{j(x)}$. Let
$$P(f)(x):=\frac1{|Q_{j(x)}|}\int_{Q_{j(x)}} f(y)dy,$$
then $P$ is linear and $P(f)$ is constant in $Q_j$, which we denote  by $q_j(f)$.
We define the linear operator $R$ by
$$R(f)=\rho_\epsilon^{\frac N2}(q_1(f),...,q_M(f))\in\R^M$$
and
$$\|v\|_2:=(\sum^M_{j=1}|v_j|^2)^{\frac12},\quad v\in \R^M.$$
We observe that for any $f\in \mathcal{T}$,
\begin{eqnarray*}
  \|P(f)\|^2_{L^2(\Omega)} &=& \sum_{j=1}^M \int_{Q_j}|P(f)(x)|^2 dx
   = \rho_\epsilon^N\sum_{j=1}^M|q_j(f)|^2
  \\&=& \|R(f)\|_2^2=\int_\Omega |f(y)|^2dy
  \\&  =& \|f\|_{L^2(\Omega)}^2\ \le C_0^2.
\end{eqnarray*}
Therefore, there exist $b_1,.....b_I\in \R^M$ such that
$$R(\mathcal{T})\subset \bigcup_{i=1}^{I}B_{\eta_\epsilon}(b_i),$$
where the balls $\{B_{\eta_\epsilon}\}$ are taken in $\R^M$.
For any $x\in\Omega$, we set
$$\beta_j(x)=\rho_\epsilon^{-\frac N2} b_{i,j(x)},$$
where $b_{i,j(x)}$ is the $j(x)$th coordinates of $b_i$.
Noticing that $\beta_j$ is constant on $Q_j$, i.e. for $x\in Q_j$, it follows that
$$P(\beta_i)(x)=\rho_\epsilon^{-\frac N2}b_{i,j}=\beta_i(x)$$
and so $q_j(\beta_i)=\rho_\epsilon^{-\frac N2}b_{i,j}$. Thus
$R(\beta_i)=b_i$.
Furthermore, for any $f\in\mathcal{T}$
\begin{eqnarray*}
 \|f-P(f)\|^2_{L^2(\Omega)} &=& \sum_{j=1}^M \int_{Q_j} |f(x)-P(f)(x)|^2 dx \\
   &=& \sum_{j=1}^M \int_{Q_j} \frac1{|Q_j|^2}|\int_{Q_j}f(x)-f(y)dy|^2 dx \\
   &\le& \frac1{\rho_\epsilon^{2N}}\sum_{j=1}^M\int_{Q_j}[\int_{Q_j}|f(x)-f(y)|dy]^2dx
\end{eqnarray*}
and for any fixed $j\in\{1,...,M\}$, by H\"{o}lder inequality, we get
\begin{eqnarray*}
&& \frac1{\rho_\epsilon^{2N}}[\int_{Q_j}|f(x)-f(y)|dy]^2 \le \frac1{\rho_\epsilon^{2N}}|Q_j|\int_{Q_j}|f(x)-f(y)|^2dy \\
   &&\qquad\qquad\le  \frac1{\rho_\epsilon^{N}}\frac1{K(2\rho_\epsilon)}\int_{Q_j}|f(x)-f(y)|^2K(x-y)dy \\
  &&\qquad\qquad\qquad\qquad\le  \frac1{\rho_\epsilon^N K(2\rho_\epsilon)}\|f\|_X^2,
\end{eqnarray*}
where $K(2\rho_\epsilon)=\inf_{|x|=2\rho_\epsilon} K(x)$.
Therefore,
\begin{eqnarray}
 \|f-P(f)\|^2_{L^2(\Omega)} \le  \frac1{\rho_\epsilon^N K(2\rho_\epsilon)} \|f\|_X^2\sum_{j=1}^M|Q_j|
   \le \frac{C}{\rho_\epsilon^N K(2\rho_\epsilon)} .\label{3.3}
\end{eqnarray}
Consequently, for any $f$, there  exists $j\in \{1,....M\}$ such that $P(f)\in B_{\eta_\epsilon}(b_j)$
and then we derive that
\begin{eqnarray*}
 &&\|f-\beta_j\|_{L^2(\Omega)}
 \\&&\quad\le \|f-P(f)\|_{L^2(\Omega)}+\|P(f)-P(\beta_j)\|_{L^2(\Omega)}+\|P(\beta_j)-\beta_j\|_{L^2(\Omega)}   \\
  \\&&\quad\le  \frac{C}{\rho_\epsilon^N K(2\rho_\epsilon)} +\frac{\|R(f)-R(\beta_j)\|_{L^2(\Omega)}}{\rho_\epsilon^{\frac N2}}\\
   \\&&\quad\le \frac{C}{\rho_\epsilon^N K(2\rho_\epsilon)}  +\frac{\eta_\epsilon}{\rho_\epsilon^{\frac N2}},
\end{eqnarray*}
where by (\ref{21-09-0}), $\frac1{(2\rho_\epsilon)^N K(2\rho_\epsilon)}\to0$ as $\epsilon\to0$ and $\frac{\eta_\epsilon}{\rho_\epsilon^{ N/2}}=\epsilon$.
As a consequence, $\mathcal{T}$ is pre-compact in $L^2(\Omega)$.

Now we are in the position to prove that $\mathcal{T}$ is pre-compact in $L^q(\Omega)$ with $q\in[1,2^*(s_0))$.
Since $L^2(\Omega)\subset L^q(\Omega)$ with $q\in[1,2)$, then $\mathcal{T}$ is pre-compact in $L^q(\Omega)$.
For $q\in(2,2^*(s_0))$, there exists $s\in(0,s_0)$ such that $q<2^*(s)$, then
using H\"{o}lder inequality with $\theta=\frac{2 (2^*(s)-q)}{q(2^*(s)-2)}$, we get that
\begin{eqnarray*}
\|f-\beta_j\|_{L^q(\Omega)}&=&\left(\int_\Omega |f-\beta_j|^{\theta q} |f-\beta_j|^{q(1-\theta)}dx \right)^{\frac1q}
 \\&\le &   \||f-\beta_j|\|^{\frac{\theta}2}_{L^2(\Omega)}\||f-\beta_j|\|^{\frac1q-\frac{\theta}2}_{L^{2^*(s)}(\Omega)}
\\&\le &\left(\frac{C}{\rho_\epsilon^N K(2\rho_\epsilon)}  +\frac{\eta_\epsilon}{\rho_\epsilon^{\frac N2}}\right)^{\frac{\theta }2}\|f\|_{X}^{\frac1q-\frac{\theta}2},
\end{eqnarray*}
thus,  $\mathcal{T}$ is pre-compact in $L^q(\Omega)$ with $q\in(2,2^*(s_0))$. The proof ends.\hfill $\Box$

\medskip

\noindent{\bf Proof of Theorem \ref{th 1}.} For Theorem \ref{th 1} part $(i)$,
let $(f_n)$ be a sequence functions in $X_0$ such that
$$\|f_n\|_X\le C,\quad \forall n\in \N$$
where $C>0$.  By Theorem \ref{sobolev}, Inequality (\ref{sobolev inequality}) follows by
(\ref{imbeding}). We obtain that the sequence $(f_n)$  is pre-compact in $L^{q}$ with $q\in[1,2^*(s_0))$, then
the compactness in Theorem \ref{th 1} follows.
\hfill $\Box$

\setcounter{equation}{0}
\section  {Existence of weak solution to (\ref{RE})}

For the proof of Theorem \ref{th 1}, we observe that
problem (\ref{RE}) has a variational structure, indeed it is
the Euler-Lagrange equation of the functional $\mathcal J:X_0\to
\R$ defined as follows
$$\mathcal J(u)=\frac 1 2 \|u\|_{X_0}^2-\frac1p\int_\Omega |u|^pdx.$$
Note  the functional $\mathcal J$ is Fr\'{e}chet differentiable in
$u\in X_0$ and for any $\varphi\in X_0$,
$$
\langle \mathcal J'(u), \varphi\rangle  = \int_Q \big(u(x)-u(y)\big)\big(\varphi(x)-\varphi(y)\big)K(x-y)dxdy
 -\int_\Omega |u|^{p-2}u(x)\varphi(x)dx.
$$

We will make use of Mountain Pass theorem to obtain the weak solution. In what follows, we check
the structure condition of Mountain Pass theorem. It is obvious that
$\mathcal J(0)=0$.
\begin{proposition}\label{pr 1}
Under the hypotheses of Theorem \ref{th 2}, there exist $\rho>0$ and
$\beta>0$ such that $\mathcal J(u)\geq \beta$, for any $u\in X_0$ with $\|u\|_{X_0}=\rho$.
\end{proposition}
{\bf Proof.}
Let $u\in X_0$, then
\begin{eqnarray*}
  \mathcal J(u) &= & \frac12\|u\|_{X_0}^2- \frac1p\int_\Omega |u(x)|^p\,dx \\
   &\ge& \frac12\|u\|_{X_0}^2- C\|u\|_{X_0}^p \\
   &=& \frac12 \|u\|_{X_0}^2(1-C\|u\|_{X_0}^{p-2}),
\end{eqnarray*}
where we used Theorem \ref{th 1} and Corollary \ref{co 1} for the inequality. We choose $\sigma>0$ such that
 $1-C\sigma^{\frac{p-2}2}=\frac12$, since $p>2$. Then for
 $\|u\|_{X_0}^2=\sigma$, $1-C\|u\|_{X_0}^{p-2}=\frac12$,
then we have
$$\mathcal J(u)\ge \frac14 \sigma.$$
The proof is complete.\hfill$\Box$

\begin{proposition}\label{pr 2}
Under the hypotheses of Theorem \ref{th 2},   there exists $e\in X_0$ such
that $\|e\|_{X_0}>\rho$ and $\mathcal
J(e)\le 0$,  where $\rho$ is given in
Proposition \ref{pr 1}.
\end{proposition}
{\bf Proof.}  We fix a function $u_0\in X_0$ with $\|u_0\|=1$ in $\Omega$. Since the space of $\{tu_0: t\in\R\}$ is
a subspace of $X_0$ with dimension 1 and all the norms are equivalent, then $\int_\Omega |u_0(x)|^p dx>0$.
Then there exists $t_0>0$ such that for $t\ge t_0$,
\begin{eqnarray*}
   \mathcal J(tu_0)&=& \frac {t^2}2\|u_0\|_{X_0}^2-\frac{t^p}p \int_\Omega |u_0(x)|^p dx\\
   &\le & C(t^2-t^p)\le  0.
\end{eqnarray*}
We choose $e=t_0u_0$.
The proof is complete.\hfill$\Box$

\medskip

We say that $\mathcal J$ has $P.S.$ condition, if for any sequence
$\{u_n\}$ in $X_0$ satisfying $\mathcal J(u_n)\to c$ and $\mathcal J'(u_n)\to0$ as $n\to\infty$, there
is a convergent subsequence, where $c\in\R$.

\begin{proposition}\label{pr 3}
Under the hypotheses of Theorem \ref{th 2},    $\mathcal J$ has  $P.S.$ condition in $X_0$.
\end{proposition}
{\bf Proof.}
Let $\{u_n\}$ be a $P.S.$ sequence, then we need to show that there are a subsequence $\{u_{n_k}\}$ and $u$ such that
$$u_{n_k}\to u\quad {\rm in}\ \ L^p(\Omega) \quad {\rm as}\ k\to\infty.$$

For some $C>0$, we have that
\begin{eqnarray}\label{4.1}
 C\|u_n\|_{X_0}\ge \mathcal J'(u_n)u_n
   =  \|u_n\|^2_{X_0} -\int_\Omega |u_n|^pdx
\end{eqnarray}
and
\begin{eqnarray}\label{4.2}
 c-1\le \mathcal J(u_n) = \frac12\|u_n\|^2_{X_0} -\frac1p\int_\Omega |u_n|^pdx.
\end{eqnarray}
Then $p\times $(\ref{4.2})-(\ref{4.1}) implies that
$$(\frac p2-1)\|u_n\|^2_{X_0}\le c+ C\|u_n\|_{X_0},$$
then $u_n$ is uniformly bounded in $X_0$.

Thus, by Theorem \ref{th 1} and Corollary \ref{co 1}, there exists a subsequence $(u_{n_k})$ and $u$ such that
$$u_{n_k}\rightharpoonup u,\quad{\rm in}\quad X_0, $$
$$u_{n_k}\to u,\quad{\rm a.e.\ in}\ \Omega\quad{\rm and\ \ in}\quad L^p(\Omega),   $$
when $k\to \infty$.
Together with $\lim_{k\to\infty} \mathcal J(u_{n_k})=c$,
we have $\|u_{n_k}\|_{X_0}\to \|u\|_{X_0}$ as $k\to\infty$.
Then we have
$u_{n_k}\to u$ in $X_0$ as $k\to\infty$.
\hfill$\Box$

\medskip

\noindent{\bf Proof of Theorem \ref{th 2}.}
By Proposition \ref{pr 1}, Proposition \ref{pr 2}
 and Proposition \ref{pr 3}, we may use Mountain Pass Theorem  (for instance, \cite[Theorem~6.1]{struwe}; see also \cite{ar, rabinowitz}) to obtain that
there exists a critical point $u\in X_0$
 of $\mathcal J$ at some value $c\ge \beta>0$. By $\beta>0$, we have $u$ is nontrivial. Therefore,
 (\ref{RE}) admits a nonnegative weak solution.
The proof is complete.\hfill$\Box$

\begin{remark}
 Suppose that $s_0\in(0,1)$ and $f:\Omega\times\R\to \R$ is a Carath\'{e}odory function
verifying the following hypothesis:
\begin{itemize}
\item[$(f_1)\ $]
\begin{enumerate}\item[]
there exist $a_1,a_2>0$ and $q\in(2,2^*(s_0))$  such that
$$|f(x,t)|\leq a_1+a_2|t|^{q-1}\quad a.e.\ x\in\Omega,\ t\in\R;$$
\end{enumerate}
\end{itemize}
\begin{itemize}
\item[$(f_2)\ $]
\begin{enumerate}\item[]
$\lim_{t\to0}\frac{f(x,t)}{|t|}=0\quad {\rm uniformly\ in}\ x\in\Omega;$
\end{enumerate}
\end{itemize}
\begin{itemize}
\item[$(f_3)\ $]
\begin{enumerate}\item[]
there exist $\mu>2$ and $r>0$ such that a.e. $x\in\Omega,t\in\R,
|t|\ge r$
$$0<\mu F(x,t)\leq tf(x,t),$$
where the function $F$ is the primitive of $f$ with respect to the variable $t$, that is
$$F(x,t)=\int_0^t f(x,\tau)d\tau.$$
\end{enumerate}
\end{itemize}
Then fractional elliptic problem (\ref{RE})
admits a nontrivial weak solution.
\end{remark}
{\bf Proof.} Using the technique in the proof of Theorem 1 in \cite{RE} and Theorem \ref{th 1} part $(ii)$,
we derive a nontrivial weak solution of (\ref{RE}) by Mountain Pass Theorem.\hfill$\Box$


\begin{thebibliography}{99}
\bibitem{ar}  A. Ambrosetti and P. Rabinowitz,
 Dual variational methods in critical point theory and applications,
{\it J. Funct. Anal. 14,}  349--381 (1973).

\bibitem{BH} A. Burchard  and H. Hajaiej,
Rearrangement inequalities for functionals with monotone integrands,
{\it J. Funct. Anal. 2,}  561-582 (2006).

\bibitem{B} H. Brezis, Functional analysis, Sobolev spaces and partial differential equations, {\it Springer} (2010).

\bibitem {CS} L. Caffarelli and L. Silvestre, An extension problem related to the fractional laplacian,
{\it Comm. Partial Differential Equations 32,}  1245-1260 (2007).

\bibitem {CS1} L. Caffarelli and  L. Silvestre, Regularity theory for fully non-linear integrodifferential equations,
{\it Comm.  Pure   Appl. Math. 62,}
597-638 (2009).

\bibitem {CS2} L. Caffarelli and L. Silvestre, Regularity results for nonlocal equations by approximation,
{\it Arch. Ration. Mech. Anal. 200(1),} 59-88 (2011).

\bibitem {CFQ} H. Chen, P. Felmer and A. Quaas, Large solution to elliptic  equations involving fractional Laplacian,
Accepted by {\it Ann. Ins.Henri Poincar\'{e}}, arXiv:1311.6044 (2013)

\bibitem {CV2} H. Chen and L. V\'{e}ron, Semilinear fractional elliptic equations  involving measures,
Accepted by {\it J. Diff. Eq.},  arXiv:1305.0945  (2013).

\bibitem {CV3} H. Chen and L. V\'{e}ron, Semilinear fractional elliptic equations with
gradient nonlinearity involving measures, {\it  J.  Funct. Anal., 266(8)}, 5467-5492 (2014).

\bibitem {CV4} H. Chen and L. V\'{e}ron, Weak and strong singular solutions of semilinear fractional  elliptic equations,
Accepted by {\it Asymp. Anal.},   arXiv:1307.7023 (2013).


\bibitem {EGE} E. Di Nazza, G.Palatucci and E. Valdinoci,
 Hitchhiker's guide to the fractional Sobolev spaces,
{\it  Bull.   Sci. Math.,   136 (5)},  521-573 (2012).

\bibitem {FQT}
P. Felmer, A. Quaas and J. Tan,
Positive solutions of non-linear Schr\"{o}dinger equation with the fractional laplacian,
{\it Proc. Roy. Soc. Edinburgh., 142}, 1237-1262 (2012).

\bibitem {FW} P. Felmer and Y. Wang, Radial symmetry of positive solutions to equations involving the fractional laplacian,
Accepted by {\it Comm. Contem. Math. }

\bibitem {HH} H. Hajaiej, Variational problems related to some fractional kinetic equations,
 arXiv:1205.1202 (2012).

\bibitem {HH1} H. Hajaiej, On the optimality of the conditions used to prove the symmetry of the minimizers
of some fractional constrained variational problems,
    {\it Ann. Inst. H. Poincar\'{e} 14(5),}  1425-1433 (2013).

\bibitem {HH2} H. Hajaiej, Existence of minimizers of functionals
involving the fractional gradient in
the abscence of compactness,
symmetry and monotonicity,
 {\it J. Math. Anal. Appl. 399(1),} 17-26 (2013).


\bibitem{rabinowitz}  P. H. Rabinowitz,
Minimax methods in critical point theory with applications to differential equations,
{\it CBMS Reg. Conf. Ser. Math., 65,  American Mathematical Society}, Providence, RI (1986).



\bibitem {RS} X. Ros-Oton and J. Serra, The Dirichlet problem for the fractional
laplacian: regularity up to the boundary,  {\it J. Math. Pures Appl., 101(3)},  275-302 (2014)

\bibitem {SVSV}
O. Savin and E. Valdinoci, Density estimates for a nonlocal variational model via the Sobolev inequality,
{\it  SIAM J. Math. Anal., 43(6)}, 2675-2687 (2011).

\bibitem {RE} R. Servadei and E. Valdinoci,
 Moutain Pass solutions for non-local elliptic operators,
 {\it J. Math. Anal. Appl. 389(2),} 887-898 (2012).

\bibitem {S1}
Y. Sire and E. Valdinoci, Fractional laplacian phase transitions and
boundary reactions: a geometric inequality and a symmetry result,
{\it J. Funct. Anal. 256}, 1842-1864 (2009).

\bibitem{struwe}  M. Struwe,
Variational methods, applications to nonlinear partial differential equations and Hamiltonian systems,
{\it Ergebnisse der Mathematik und ihrer Grenzgebiete,  3,  Springer Verlag}, Berlin--Heidelberg (1990).




\end{thebibliography}
\end{document}